\documentclass[a4paper,reqno]{amsart}
\usepackage{amssymb, amsmath, amscd}
\usepackage[final]{graphicx}
\usepackage{float}\usepackage{wrapfig}
\usepackage{color}
\usepackage{bbm}
\usepackage[final]{graphicx}

\def\Tr{\operatorname{Tr}}
\newtheorem{theorem}{Theorem}[section]

\newtheorem{lemma}[theorem]{Lemma}
\newtheorem{remark}[theorem]{Remark}

\makeatletter
 \@addtoreset{equation}{section}
\makeatother
\begin{document}
\title
{Harmonic spaces}
\author
{P. B. Gilkey and J. H. Park}
\address{PG: Mathematics Department, University of Oregon, Eugene OR 97403-1222, USA}
\email{gilkey@uoregon.edu}
\address{JHP: Department of
Mathematics, Sungkyunkwan University, Suwon, 16419 Korea.}
\email{parkj@skku.edu}
\begin{abstract}We use the density function of a harmonic space to obtain estimates for the eigenvalues of the Jacobi operator; when
these estimates are sharp, then the harmonic space is a symmetric Osserman space.\end{abstract}
\maketitle

\section{Introduction}
\subsection{The density function} If $\mathcal{M}=(M,g)$ is a
Riemannian manifold of dimension $m$, let $\operatorname{dvol}_{\mathcal{M}}$ be the associated measure.
We identify a neighborhood of $P$ with a neighborhood of $0$ in $T_PM$ using the exponential map
$\exp_P:T_PM\rightarrow M$.
Since the tangent space $T_PM$ is an inner product space,
it inherits a natural measure we denote by $d\xi$. We fix a local orthonormal frame for $T_PM$ to introduce
coordinates and express $d\xi=d\xi^1\dots d\xi^m$. We may then express
$$
\exp_P^*\{\operatorname{dvol}_{\mathcal{M}}\}(\xi)=\tilde\Theta_{\mathcal{M}}(P;\xi)d\xi\text{ for }\tilde\Theta=\sqrt{\det(g_{ij})}\,.
$$
In such a coordinate system, $g_{ij}=\delta_{ij}+O(\xi^2)$ and thus
$\tilde\Theta_{\mathcal{M}}(P;\xi)=1+O(\xi^2)$. We expand $\tilde\Theta$ in a formal Taylor series
$$
\tilde\Theta_{\mathcal{M}}(P;\xi)\sim 1+\sum_{k=2}^\infty\mathcal{H}_k(P;\xi)
$$
where $\mathcal{H}_k(\xi)$ is homogeneous of degree $k$ in $\xi$, i.e.
one has that $\mathcal{H}_k(c\xi)=c^k\mathcal{H}_k(\xi)$.
To simplify the notation, let
$\mathcal{J}_k(\xi)$ be the endomorphism of $T_PM$ defined by the relationship:
$$
g(\mathcal{J}_k(\xi)\xi_1,\xi_2)=R(\xi_1,\xi,\xi,\xi_2;\xi\dots\xi)\,.
$$
In brief, $\mathcal{J}=\mathcal{J}_0$ is the Jacobi operator and $\mathcal{J}_k(\xi)=\nabla_\xi^k\mathcal{J}(\xi)$.
We will establish the following result in Section~\ref{S2}.
\begin{theorem}\label{T1.1} Let $\mathcal{M}$ be a Riemannian manifold. Then
\begin{enumerate}
\smallbreak\item $\displaystyle\mathcal{H}_2(\xi)=-\frac{\Tr\{\mathcal{J}(\xi)\}}{6}$.
\smallbreak\item $\displaystyle\mathcal{H}_3(\xi)=-\frac{\Tr\{\mathcal{J}_1(\xi)\}}{12}$.
\smallbreak\item $\displaystyle\mathcal{H}_4(\xi)=\frac{\Tr\{\mathcal{J}(\xi)\}^2}{72}
-\frac{\Tr\{\mathcal{J}(\xi)^2\}}{180}-\frac{\Tr\{\mathcal{J}_2(\xi)\}}{40}$
 \smallbreak\item $\displaystyle\mathcal{H}_5(\xi)=\frac{\Tr\{\mathcal{J}(\xi)\} \Tr\{\mathcal{J}_1(\xi)\}}{72}
 -\frac{\Tr\{\mathcal{J}(\xi)\mathcal{J}_1(\xi)\}}{180}-\frac{\Tr\{\mathcal{J}_3(\xi)\}}{180}$.
 \smallbreak\item $\displaystyle\mathcal{H}_6(\xi)=-\frac{\Tr\{\mathcal{J}(\xi)\}^3}{1296}
 +\frac{\Tr\{\mathcal{J}(\xi)\} \Tr\{\mathcal{J}(\xi)^2\}}{1080}+\frac{\Tr\{\mathcal{J}(\xi)\} \Tr\{\mathcal{J}_2(\xi)\}}{240}$
 \smallbreak\qquad$\displaystyle
 -\frac{\Tr\{\mathcal{J}(\xi)^3\}}{2835}
 -\frac{\Tr\{\mathcal{J}(\xi)\mathcal{J}_2(\xi)\}}{630}+\frac{\Tr\{\mathcal{J}_1(\xi)\}^2}{288}
 -\frac{\Tr\{\mathcal{J}_1(\xi)^2\}}{672}$
 \smallbreak\qquad$\displaystyle-\frac{\Tr\{\mathcal{J}_4(\xi)\}}{1008}$.
 \goodbreak\smallbreak\item $\displaystyle\mathcal{H}_7(\xi)=\frac{\Tr\{\mathcal{J}(\xi)\} \Tr\{\mathcal{J}(\xi)\mathcal{J}_1(\xi)\}}{1080}
 -\frac{\Tr\{\mathcal{J}(\xi)\}^2 \Tr\{\mathcal{J}_1(\xi)\}}{864}-\frac{\Tr\{\mathcal{J}_5(\xi)\}}{6720}$
 \smallbreak\qquad$\displaystyle+\frac{\Tr\{\mathcal{J}(\xi)\} \Tr\{\mathcal{J}_3(\xi)\}}{1080}
 +\frac{\Tr\{\mathcal{J}(\xi)^2\} \Tr\{\mathcal{J}_1(\xi)\}}{2160}-\frac{\Tr\{\mathcal{J}(\xi)^2\mathcal{J}_1(\xi)\}}{1890}$
 \smallbreak\qquad$\displaystyle
 -\frac{\Tr\{\mathcal{J}(\xi)\mathcal{J}_3(\xi)\}}{3024}+\frac{\Tr\{\mathcal{J}_1(\xi)\} \Tr\{\mathcal{J}_2(\xi)\}}{480}
 -\frac{\Tr\{\mathcal{J}_1(\xi)\mathcal{J}_2(\xi)\}}{1120}$.
 \smallbreak\item $\displaystyle\mathcal{H}_8(\xi)=\frac{\Tr\{\mathcal{J}(\xi)\}^4}{31104}-\frac{\Tr\{\mathcal{J}(\xi)\}^2 \Tr\{\mathcal{J}(\xi)^2\}}{12960}
-\frac{\Tr\{\mathcal{J}(\xi)\}^2 \Tr\{\mathcal{J}_2(\xi)\}}{2880}$
\smallbreak\qquad$\displaystyle+\frac{\Tr\{\mathcal{J}(\xi)\} \Tr\{\mathcal{J}(\xi)^3\}}{17010}
+\frac{\Tr\{\mathcal{J}(\xi)\} \Tr\{\mathcal{J}(\xi)\mathcal{J}_2(\xi)\}}{3780}-\frac{\Tr\{\mathcal{J}_6(\xi)\}}{51840}$
\smallbreak\qquad$\displaystyle
-\frac{\Tr\{\mathcal{J}(\xi)\} \Tr\{\mathcal{J}_1(\xi)\}^2}{1728}
+\frac{\Tr\{\mathcal{J}(\xi)\} \Tr\{\mathcal{J}_1(\xi)^2\}}{4032}$
\smallbreak\qquad$\displaystyle-\frac{\Tr\{\mathcal{J}_2(\xi)^2\}}{7200}
+\frac{\Tr\{\mathcal{J}(\xi)\} \Tr\{\mathcal{J}_4(\xi)\}}{6048}
+\frac{\Tr\{\mathcal{J}(\xi)^2\} \Tr\{\mathcal{J}_2(\xi)\}}{7200}$
\smallbreak\qquad$\displaystyle+\frac{\Tr\{\mathcal{J}(\xi)^2\} ^2}{64800}
-\frac{\Tr\{\mathcal{J}(\xi)^4\}}{37800}-\frac{17 \Tr\{\mathcal{J}(\xi)^2\mathcal{J}_2(\xi)\}}{113400}$
\smallbreak\qquad$\displaystyle
+\frac{\Tr\{\mathcal{J}(\xi)\mathcal{J}_1(\xi)\} \Tr\{\mathcal{J}_1(\xi)\}}{2160}-\frac{5 \Tr\{\mathcal{J}(\xi)\mathcal{J}_1(\xi)^2\}}{18144}
-\frac{\Tr\{\mathcal{J}(\xi)\mathcal{J}_4(\xi)\}}{18144}$
\smallbreak\qquad$\displaystyle+\frac{\Tr\{\mathcal{J}_1(\xi)\} \Tr\{\mathcal{J}_3(\xi)\}}{2160}
-\frac{\Tr\{\mathcal{J}_1(\xi)\mathcal{J}_3(\xi)\}}{5184}+\frac{\Tr\{\mathcal{J}_2(\xi)\}^2}{3200}$.
 \end{enumerate}\end{theorem}
 \begin{remark}\rm Our formulas agree with those in Gray~\cite{Gr74} for $2\le n\le 6$;
the formulas for $\mathcal{H}_7$ and $\mathcal{H}_8$
are, we believe, new. The following formula of Vanhecke~\cite{V} can be used to determine $\mathcal{H}_9$
from these formulas; the odd asymptotic coefficients are expressible in terms of the even coefficients:
$$
\mathcal{H}_{k}(P;\xi)=\frac12\sum_{i=2}^{k-1}\frac{(-1)^i}{(k-i)!}
\left(\nabla_{\xi\dots\xi}^{k-i}\mathcal{H}_i\right)(P;\xi)\text{ for }k\text{ odd}\,.
$$
\end{remark}

We can introduce geodesic polar coordinates identifying $\xi=r\cdot\theta$ where $r>0$ and $\|\theta\|=1$.  We then have
$d\xi=r^{m-1}drd\theta$ and $\exp_P^*\operatorname{dvol}_{\mathcal{M}}=\Theta_{\mathcal{M}}(P;\xi)r^{m-1}drd\theta$ where
$\Theta_{\mathcal{M}}:=r^{1-m}\tilde\Theta_{\mathcal{M}}$. We shall work with $\Theta$ instead of $\tilde\Theta$ hence forth
for historical reasons, but both $\Theta$ and $\tilde\Theta$ contain the same information.

\subsection{Harmonic spaces}
Let $\Omega_{\mathcal{M}}(P,Q):=\frac12 d_{\mathcal{M}}(P,Q)^2$ where $d_{\mathcal{M}}(P,Q)$ is the geodesic
distance from $P$ to $Q$.  Let $S_{\mathcal{M}}(P;r)=\{Q:d_{\mathcal{M}}(P,Q)=r\}$ be the geodesic sphere of radius $r$ centered at $P$.
If $f$ is a smooth function, let $M_{\mathcal{M}}(f;P,r)$
be the average value of a function $f$ on the geodesic sphere of radius $r$ centered about $P$:
$$
M_{\mathcal{M}}(f;P,r):=\frac1{\operatorname{vol}(S_{\mathcal{M}}(P;r))}\int_{S_{\mathcal{M}}(P;r)}f\operatorname{dvol}_{\mathcal{M}}
$$
We have the following result, see Berndt, Tricerri, and Vanhecke~\cite{BTV78} and Besse~\cite{B78}.

\begin{theorem}
The following conditions are equivalent and if any is satisfied, then $\mathcal{M}$ is said to be a harmonic space.
\begin{enumerate}
\item $\Theta_{\mathcal{M}}(P;\xi)=\Theta_{\mathcal{M}}(P;\|\xi\|)$,
i.e. $\Theta_{\mathcal{M}}$ only depends on $P$ and $\|\xi\|$.
\item $\Theta_{\mathcal{M}}(P;\xi)=\Theta_{\mathcal{M}}(\|\xi\|)$, i.e. $\Theta_{\mathcal{M}}$ only depends on $\|\xi\|$.
\item Near any $P\in M$, there exists a non-constant harmonic radial function.
\item Near each point $P$ of $M$, $\Delta\Omega_{\mathcal{M}}(P,\cdot)$ is a function of $\Omega_{\mathcal{M}}(P,\cdot)$.
\item Small geodesic spheres in $M$ have constant mean curvature.
\item If $f$ is harmonic near $P$, then $M_{\mathcal{M}}(f;P,r)=f(P)$ for all $P$ in $M$.
\end{enumerate}
\end{theorem}

If we assume the $\mathcal{M}$ is a harmonic space, then the invariants
$\mathcal{H}_k$ are constant on the unit sphere bundle. This implies, in particular, by Theorem~\ref{T1.1} that
$\mathcal{M}$ is 2-stein, i.e. $\Tr\{\mathcal{J}(\xi)\}$ and $\Tr\{\mathcal{J}(\xi)^2\}$ are constant on the unit tangent
bundle. Differentiating these relations yields additional relations and permits us to simplify the formulas of Theorem~\ref{T1.1}
in this situation. We will establish the following result in Section~\ref{S3}.
\begin{theorem}\label{T1.4}
 Let $\mathcal{M}$ be a harmonic space. Then $\mathcal{H}_{2k+1}(\xi)=0$ for all $k$.
\begin{enumerate}
\item $\displaystyle\mathcal{H}_2(\xi)=-\frac16\Tr\{\mathcal{J}(\xi)\}$.
\smallbreak\item $\displaystyle\mathcal{H}_4(\xi)=\frac1{72}\Tr\{\mathcal{J}(\xi)\}^2-\frac1{180}\Tr\{\mathcal{J}(\xi)^2\}$.
\smallbreak\item $\displaystyle\mathcal{H}_6(\xi)=-\frac{\Tr\{\mathcal{J}(\xi)\}^3}{1296}
+\frac{\Tr\{\mathcal{J}(\xi)\}\Tr\{\mathcal{J}(\xi)\}^2}{1080}
-\frac{\Tr\{\mathcal{J}(\xi)^3\}}{2835}$
$\displaystyle+\frac{\Tr\{\mathcal{J}_1(\xi)^2\}}{10080}$.
\smallbreak\item $\displaystyle\mathcal{H}_8(\xi)=\frac{\Tr\{\mathcal{J}(\xi)\}^4}{31104}-\frac{\Tr\{\mathcal{J}(\xi)\}^2 \Tr\{\mathcal{J}(\xi)^2\}}{12960}+\frac{\Tr\{\mathcal{J}(\xi)\} \Tr\{\mathcal{J}(\xi)^3\}}{17010}$
\smallbreak\qquad$\displaystyle-\frac{\Tr\{\mathcal{J}(\xi)\} \Tr\{\mathcal{J}_1(\xi)^2\}}{60480}
+\frac{\Tr\{\mathcal{J}(\xi)^2\}^2}{64800}-\frac{\Tr\{\mathcal{J}(\xi)^4\}}{37800}$
\smallbreak\qquad$\displaystyle-\frac{ \Tr\{\mathcal{J}(\xi)^2\mathcal{J}_2(\xi)\}}{340200}+\frac{\Tr\{\mathcal{J}(\xi)\mathcal{J}_1(\xi)^2\}}{54432}-\frac{\Tr\{\mathcal{J}_2(\xi)^2\}}{907200}$.
\end{enumerate}
\end{theorem}

\begin{remark}\rm
Assertions~(1), (2), and (3) are not new; they follow from the work of
Besse~\cite{B78}, Copson and Ruse~\cite{CR40}, Gray~\cite{Gr74}, Ledger~\cite{L57}, and Lichnerowicz~\cite{L44}.
Assertion~(4) is, we believe, new.
\end{remark}

As noted above, any harmonic space is Einstein.
It follows from work of Ramachandran and Ranjan \cite{RR97} and the work of
Ranjan and Shah \cite{RS2002} that any Ricci flat harmonic space is flat;
thus we may assume that the scalar curvature $\tau$ is non-zero. This naturally divides
the discussion into the case of constant positive scalar curvature and constant negative scalar curvature.
We can rescale the metric to let $\mathcal{M}(c):=(M,c^2g)$ for $c>0$.
If $\mathcal{M}$ is harmonic, then $\mathcal{M}(c)$ is harmonic and
\begin{equation}\label{Eq1.a}
\Theta_{\mathcal{M}(c)}(r)=c^{1-m}\Theta_{\mathcal{M}}(cr)\,.
\end{equation}
The following density functions will play a central role in what follows. Set
\begin{equation}\label{E1.b}
\Theta_\varepsilon(m,k)(r):=\left\{\begin{array}{cl}
\sin(r)^{m-1}\cos(r)^k&\text{if }\varepsilon=+\\
\sinh(r)^{m-1}\cosh(r)^k&\text{if }\varepsilon=-\end{array}\right\}\,.
\end{equation}
Apart from the flat case where $\Theta_{\mathcal{M}}(r)=r^{m-1}$, these functions are
the only known density functions for a harmonic space modulo the rescaling discussed
in Equation~(\ref{Eq1.a}).

\subsection{Two point homogeneous spaces}
A Riemannian manifold
$\mathcal{M}$ is said to be a {\it two point homogeneous space} if the
group of isometries acts transitively on the associated the unit tangent bundle. In this instance,
$\Theta_{\mathcal{M}}(P;\xi)=\Theta_{\mathcal{M}}(\|\xi\|)$
so any two point homogeneous space is a harmonic space. The two point homogeneous spaces
have been classified; such a manifold is either flat or is a rank one symmetric space.

Let $R_{\mathcal{M}}(\xi,\eta):=\nabla_\xi\nabla_\eta-\nabla_\eta\nabla_\xi-\nabla_{[\xi,\eta]}$ be the associated curvature operator
and let $R_{\mathcal{M}}(\xi_1,\xi_2,\xi_3,\xi_4):=g(R_{\mathcal{M}}(\xi_1,\xi_2)\xi_3,\xi_4)$ be the associated curvature tensor.
One defines the {\it Jacobi operator} $\mathcal{J}_{\mathcal{M}}$ by setting $\mathcal{J}_{\mathcal{M}}(\xi)\eta:=R_{\mathcal{M}}(\eta,\xi)\xi$;
this is a self-adjoint
operator on the tangent bundle of $M$. The Ricci tensor $\rho_{\mathcal{M}}$ satisfies
$\rho_{\mathcal{M}}(\xi,\xi)=\operatorname{Tr}(\mathcal{J}_{\mathcal{M}}(\xi))$. Since $g(R_{\mathcal{M}}(\eta,\xi)\xi,\xi)=0$,
$\mathcal{J}_{\mathcal{M}}(\xi)$ preserves $\xi^\perp$;
we define the {\it reduced Jacobi operator} $\tilde{\mathcal{J}_{\mathcal{M}}}(\xi)$
to be the restriction of $\mathcal{J}_{\mathcal{M}}(\xi)$ to $\xi^\perp$;
the eigenvalues of $\mathcal{J}_{\mathcal{M}}(\xi)$, counted with
multiplicity, are the eigenvalues of $\tilde{\mathcal{J}_{\mathcal{M}}}(\xi)$ plus an additional ``trivial" eigenvalue $0$.

If $\mathcal{M}$ is a two point homogeneous space, then the eigenvalues of the Jacobi operator are constant
on the unit tangent bundle. Osserman~\cite{oss} wondered if the converse held; if the eigenvalues of the Jacobi
operator are constant on the unit tangent bundle, is $\mathcal{M}$ locally isometric to a two point homogeneous space;
such manifolds have been called {\it Osserman manifolds} by later authors. This question by Osserman has been
settled in the affirmative by Chi~\cite{C88} and Nikolayevsky~\cite{niko3,niko4}
 if $m\ne16$; the case $m=16$ is still open and this plays a role in our analysis.

Let $\mathcal{M}$ be a two point homogeneous space. We summarize the information we shall need about
the density function $\Theta_{\mathcal{M}}$ and about
the eigenvalues of the Jacobi operator as follows.
Let $\xi$ be a unit tangent vector. We may decompose $\xi^\perp=V_1\oplus V_2$
where $\dim(V_1)=\nu_1$ and $\dim(V_2)=\nu_2$ are constant.
Then $J(\xi)\eta=\kappa_i\eta$ if $\eta\in V_i$, i.e. $\mathcal{M}$ is an {\it Osserman space} and the reduced Jacobi operator has
at most two distinct eigenvalues $\kappa_1$ and $\kappa_2$ with multiplicities $\nu_1$ and $\nu_2$, respectively.
In positive curvature, up to the rescaling $\mathcal{M}\rightarrow\mathcal{M}(c)$, the rank one  symmetric spaces are the sphere $S^n$, complex
projective space $\mathbb{CP}^n$, quaternonic projective space $\mathbb{HP}^n$, and the {Cayley plane (or octonionic projective plane)} $\mathbb{OP}^2$.
There are corresponding negative curvature duals: hyperbolic space $\mathbb{H}^n$ is the
negative curvature of the sphere, $\widetilde{\mathbb{CP}^n}$ is the negative curvature dual of complex projective space,
$\widetilde{\mathbb{HP}^n}$ is the negative curvature dual of quaternionic projective space, and
$\widetilde{\mathbb{OP}^2}$ is the negative curvature dual of the Cayley plane.

We refer to Watanabe~\cite{W75,W76}
for the calculation of the characteristic function of the spaces other than the Cayley plane and its non-compact
dual and to Euh, Park, and Sekigawa~\cite{EPS} for the corresponding calculation for the Cayley plane and its
non-compact dual;
the curvature tensors are well known classically. We have $\Theta_{\mathcal{M}}=\Theta_\varepsilon(m,k)$ where
$k=\nu_1$.
\smallbreak\centerline{Positive scalar curvature}
$$\begin{array}{| l | c | r | r | r | r | c |}\noalign{\hrule}
&\dim&k=\nu_1&\kappa_1&\nu_2&\kappa_2&\varepsilon\\\noalign{\hrule}
S^n&m=n&k=0&&n-1&1&\varepsilon=+\\\noalign{\hrule}
\mathbb{CP}^n&m=2n&k=1&4&2n-2&1&\varepsilon=+\\ \noalign{\hrule}
\mathbb{HP}^n&m=4n&k=3&4&4n-4&1&\varepsilon=+\\ \noalign{\hrule}
\mathbb{OP}^2&m=16&k=7&4&8&1&\varepsilon=+\\ \noalign{\hrule}
\end{array}$$
\centerline{Negative scalar curvature}
$$\begin{array}{| l | c | r | r | r | r | c |}\noalign{\hrule}
&\dim&k=\nu_1&\kappa_1&\nu_2&\kappa_2&\varepsilon\\\noalign{\hrule}
\widetilde{S^n}=\mathbb{H}^n&m=n&k=0&&n-1&-1&\varepsilon=-\\\noalign{\hrule}\vphantom{\vrule height 12pt}
\widetilde{\mathbb{CP}^n}&m=2n&k=1&-4&2n-2&-1&\varepsilon=-\\ \noalign{\hrule}\vphantom{\vrule height 12pt}
\widetilde{\mathbb{HP}^n}&m=4n&k=3&-4&4n-4&-1&\varepsilon=-\\ \noalign{\hrule}\vphantom{\vrule height 12pt}
\widetilde{\mathbb{OP}^2}&m=16&k=7&-4&8&-1&\varepsilon=-\\ \noalign{\hrule}\end{array}$$

\subsection{Damek-Ricci spaces}
Let $\mathfrak{n}=\mathfrak{z}\oplus\mathfrak{u}$ be an orthogonal decomposition of a
nilpotent Lie algebra $\mathfrak{n}$ with center $\mathfrak{z}$ such that
$[\mathfrak{u},\mathfrak{u}]\subset\mathfrak{z}$, i.e. $\mathfrak{n}$ is 2-step nilpotent.
Let $\langle\cdot,\rangle\cdot$ be a left invariant metric on $\mathfrak{n}$.
For $z\in\mathfrak{z}$, define a skew-symmetric operator $J_z:\mathfrak{u}\rightarrow\mathfrak{u}$
by $\langle J_zx,y\rangle=\langle z,[x,y]\rangle$ for $x,y\in\mathfrak{u}$. Let $\mathfrak{g}=\mathfrak{n}\oplus\mathfrak{a}$ where
$\mathfrak{a}$ is 1-dimensional and orthogonal to $\mathfrak{n}$ spanned by a unit
vector $A$ where
$\operatorname{ad}_{A|\mathfrak{u}}=\frac12\operatorname{id}_{\mathfrak{u}}$, and
$\operatorname{ad}_{A|\mathfrak{z}}=\operatorname{id}_{\mathfrak{z}}$ and where
$J_z^2=-\|z\|^2\operatorname{id}_{\mathfrak{u}}$ for all $z\in\mathfrak{z}$.

Following the seminal paper of Damek and Ricci~\cite{AD},
the associated Lie group, which is a {\it solvmanifold},  is called a {\it Damek--Ricci space}. Such a manifold
is a harmonic space of dimension $m=\dim(\mathfrak{z})+\dim(\mathfrak{u})+1$; with a suitably normalized metric,
the volume element takes the form
$$
\Theta_{\mathcal{M}}=\Theta_-(m=\dim(\mathfrak{z})+\dim(\mathfrak{u})+1,k=\dim(\mathfrak{z}))\,.
$$
Rank one non-compact symmetric spaces (including the real hyperbolic space, if one allows $\mathfrak{u} = 0$)
are specific cases of Damek-Ricci spaces. If $\mathcal{M}$ is a rank one  symmetric space of non-compact type, then
$$\begin{array}{| c | c | c | c | c | c |}\noalign{\hrule}
\mathcal{M}&\dim(\mathcal{M})&\dim\{\mathfrak{z}\}&\mathcal{M}&
\dim(\mathcal{M})&\dim\{\mathfrak{z}\}\\\noalign{\hrule}
\mathbb{H}^m&m&0&\widetilde{\mathbb{CP}}^n&2n&1\\\noalign{\hrule}
\widetilde{\mathbb{HP}}^n&4n&3&\widetilde{\mathbb{OP}}^{2}&16&7\\\noalign{\hrule}
\end{array}$$
There exist examples of Damek--Ricci spaces with $\dim(\mathfrak{z})$ arbitrary. The first non-symmetric Damek-Ricci space is
7-dimensional; for $\dim(\mathfrak{z})$ small, the dimensions of the corresponding non-symmetric Damek--Ricci spaces are
($n=0,1,\dots$)
$$\begin{array}{| c | c | c | c | c | c | c | c |}\noalign{\hrule}
\dim(\mathfrak{z})&\dim(\mathcal{M})&\dim(\mathfrak{z})&\dim(\mathcal{M})
&\dim(\mathfrak{z})&\dim(\mathcal{M})&\dim(\mathfrak{z})&\dim(\mathcal{M})
\\ \noalign{\hrule}
1&&2&7+4n&3&12+4n&4&13+8n\\ \noalign{\hrule}
5&14+8n&6&15+8n&7&24+8n&8&25+16n\\\noalign{\hrule}\end{array}$$

Note that no non-symmetric Damek--Ricci space can be Osserman (see Theorem~2 of \cite{BTV78});
in particular, the eigenvalues of the Jacobi operator of a non-symmetric Damek--Ricci space are not constant on the unit tangent bundle.

\subsection{Classification results}
In dimensions $2$ and $3$ every harmonic space is a space form, i.e. has constant sectional curvature.
Lichnerowicz~\cite{L44} conjectured that a 4-dimensional harmonic space was a locally symmetric space; this
was established by Walker~\cite{W49}. Ledger~\cite{L57} conjectured that every harmonic space is either
locally flat or locally a rank one  symmetric space. Szab\'o~\cite{S90} showed this was true if $\mathcal{M}$
was complete and had positive scalar curvature; since any harmonic space is Einstein, this implies $\mathcal{M}$
is compact by Myer's Theorem.
Besson et al.~\cite{BG95} established the conjecture if $\mathcal{M}$ is compact with negative sectional
curvature. Nikolayevsky~\cite{N05} showed every harmonic space of dimension 5 is a space form.
However, the conjecture of Lichnerowicz is false in the negative curvature setting without the
hypothesis of compactness. As noted above, Damek and Ricci~\cite{DR57}
constructed complete, non-compact harmonic spaces which are not rank
one symmetric spaces for any dimension  $m\geq 7$; the conjecture is still open in dimension $6$.

Heber~\cite{H06} showed that any simply connected homogeneous harmonic space is flat, or rank-one symmetric,
or a non-symmetric Damek-Ricci space.
It is not known if there are inhomogeneous examples.
It is an open question as to whether there exist germs of
harmonic spaces which are essentially geodesically incomplete, i.e. which do not embed in
complete harmonic space; the research of the present paper is motivated, at least in part, by
this question. We also refer to work of Csik\`os and Horv\'ath~\cite{CH16,CH18} characterizing
harmonic manifolds in terms of the total mean curvature and the total scalar curvature of tubular
hypersurfaces about curves.

\subsection{Characterization of harmonic spaces by their density function}
The following result by Ramachandran and Ranjan~\cite{RR97} is evocative.
\begin{theorem}\label{T1.6}
Let $\mathcal{M}$ be a non-compact simply connected complete harmonic space of dimension $m$
with $\Theta_{\mathcal{M}}=\Theta_-(m,k)$.
\begin{enumerate}
\item If $k=0$, then $\mathcal{M}$ is isometric to $\mathbb{H}^m$.
\item If $k=1$ and if $\mathcal{M}$ is K\"ahler, then $m=2\bar m$ and
$\mathcal{M}$ is isometric to $\widetilde{\mathbb{CP}}^{\bar m}$.
\item If $k=3$ and if $\mathcal{M}$ is quaternionic K\"ahler, then
$m=4\bar m$ and $\mathcal{M}$ is isometric to $\widetilde{\mathbb{HP}^{\bar m}}$.
\end{enumerate}
\end{theorem}

\subsection{The density functions and the eigenvalues of the Jacobi operator}
We say that a Riemannian manifold $\mathcal{M}$ is $k$-stein if $\operatorname{Tr}(\mathcal{J}_{\mathcal{M}}(\xi)^\ell)$ is constant
on the unit tangent bundle for $1\le\ell\le k$. The $1$-stein condition is the familiar Einstein condition, $\rho_{\mathcal{M}}(\xi,\xi)=c\|\xi\|^2$;
polarizing this identity shows $\rho_{\mathcal{M}}=cg$.
If $A$ and $B$ are two symmetric $\ell\times\ell$ matrices, then the eigenvalues, counted with multiplicity,
of $A$ and $B$ coincide if and only if $\operatorname{Tr}(A^i)=\operatorname{Tr}(B^i)$ for $1\le i\le\ell$. Since the trivial
eigenvalue $0$ plays no role, a manifold of dimension $m$ is Osserman if and only if it is $(m-1)$-stein.

Nikolaveysky~\cite{N05} showed that the volume density function of a harmonic manifold is an
exponential polynomial, i.e. a finite linear combination of terms of the form $\Re(c_ie^{\lambda_ir}r^{m_i})$.
Motivated by this result, we shall assume henceforth that $\Theta_{\mathcal{M}}=\Theta_\varepsilon(m,k)$
as these are the only known examples. We observe that in the indefinite signature
setting, the volume density function can have the form $\sinh(r)^{p-1}\sin(r)^{q}\cosh(r)^k$ so the situation
is more complicated there.

We refer to Berndt, Tricerri, and Vanhecke~\cite{BTV78} for a computation of the eigenvalues of the
Jacobi operator for an arbitrary Damek--Ricci space. We are working in a quite general context and
do not assume $\mathcal{M}$ is a Damek-Ricci space.
We will use Theorem~\ref{T1.4} to establish the following result in Section~\ref{S4.1}; {Assertion~3} generalizes
Theorem~\ref{T1.6}~(1) to the positive curvature setting.

\begin{lemma}\label{L1.7}
Let $\mathcal{M}$ be a harmonic space with $\Theta_{\mathcal{M}}=\Theta_\varepsilon(m,k)$.
\begin{enumerate}
\item $m=\dim(\mathcal{M})$ and $0\le k\le m-1$.
\item If $k=0$, then $\mathcal{M}$ has constant
sectional curvature $\varepsilon$.
\item If $k=m-1$, then $\mathcal{M}$ has constant
sectional curvature $4\varepsilon$.
\end{enumerate}\end{lemma}

By Lemma~\ref{L1.7}, if $\Theta_{\mathcal{M}}=\Theta_\varepsilon(m,k)$, then
we may assume that $1\le k<m-1$.
We will establish the following result in Section~\ref{S4.2}.

\begin{theorem}\label{T1.8}
Let  $1\le k<m-1$.
Let $\mathcal{M}$ be a harmonic space. Assume that $\Theta_{\mathcal{M}}=\Theta_\varepsilon(m,k)$.
Let $\xi$ be a unit tangent vector in  $\mathcal{M}$. Order the eigenvalues of the reduced Jacobi operator
$\tilde{\mathcal{J}_{\mathcal{M}}}(\xi)$
so $\varepsilon\lambda_1\ge\varepsilon\lambda_2\ge\dots\ge\varepsilon\lambda_{m-1}$.
\begin{enumerate}
\item We have $\varepsilon\lambda_1+\dots+\varepsilon\lambda_k\le4k$ and
$\varepsilon\lambda_{k+1}+\dots+\varepsilon\lambda_{m-1}\ge m-1-k$.
\item If either of the inequalities in Assertion~1 is an equality for all unit vectors $\xi$,
then $\mathcal{M}$ is a rank one symmetric space.
\end{enumerate}\end{theorem}

Note that Theorem~\ref{T1.8} is a purely local result; we impose no assumptions of completeness
or other geometric hypotheses.
This gives an abstract characterization of the harmonic spaces which are rank 1 symmetric spaces
since the inequalities of Assertion~1 are in fact equalities in this instance.

Theorem~\ref{T1.8} result in gives an estimate on averages of the eigenvalues of the reduced Jacobi operator.
One can also estimate
individual eigenvalues.

\begin{theorem}\label{T1.9}
Let $\Theta_{\mathcal{M}}=\Theta_\varepsilon(m,k)$. Let $\lambda$ be an eigenvalue of the reduced Jacobi operator.
Assume $m>2$. Then
\begin{eqnarray*}
&&\frac{m+3k-1}{m-1}-\mathcal{E}(m,k)\le\varepsilon\lambda\le\frac{m+3k-1}{m-1}+\mathcal{E}(m,k)\quad\text{for}\\
&&\mathcal{E}(m,k):=\frac{3m-6}{m-1}\sqrt{\frac{k(m-1-k)}{m-2}}\,.
\end{eqnarray*}
\end{theorem}

\section{The asymptotic series of the density function}\label{S2}
This section is devoted to the proof of Theorem~\ref{T1.1}.
Let $\mathcal{M}$ be an arbitrary Riemannian manifold. Fix a point $P$ of $M$ and let $\xi\in T_PM$. We use
geodesic polar coordinates relative to some local orthonormal frame for $T_PM$ and expand
$$
\Theta_P(\xi)\sim r^{m-1}\left(1+\sum_{k=2}^\infty\mathcal{H}_k(\xi)\right)\quad\text{where}\quad r=\|\xi\|\,.
$$
We have that
$\mathcal{H}_k(\xi)$ is homogeneous of degree $k$ in the derivatives of the metric evaluated at $P$ and that
$\mathcal{H}_k(r\xi)=r^k\mathcal{H}_k(\xi)$. Since the derivatives of the metric can be expressed in terms of curvature at
the origin, we have $\mathcal{H}_k(\xi)$ is expressible in terms of curvature.
Let $I=(i_1,\dots,i_\ell)$ be sequence of non-negative integers listed in increasing order. Set
$$
\mathcal{J}_I\xi)=\Tr\{\mathcal{J}_{i_1}(\xi)\dots\mathcal{J}_{i_\ell}(\xi)\}\text{ and }
\operatorname{ord}(I)=(2+i_1)+\dots+(2+i_\ell)\,.
$$
We then have that $\mathcal{J}_I$ is homogeneous of degree $\operatorname{ord}(I)$ in the derivatives of the metric.
We use H. Weyl's theorem on the invariants of the orthogonal group to see that
$\mathcal{H}_k$ is a polynomial of total order $k$ in the variables $\mathcal{J}_I$.
A-priori, the coefficients in this polynomial depend on the dimension. Standard arguments using taking
products with the circle enable one to dimension shift and show that the coefficients can actually be taken
to be dimension independent; we refer, for example, to the discussion in Gilkey~\cite{G95} concerning heat
trace asymptotics. We adopt the {\it Einstein convention} and sum over repeated indices relative to a local orthonormal frame
$\{e_1,\dots,e_m\}$ for the tangent bundle to expand
\begin{equation}\label{E2.a}\begin{array}{lll}
\mathcal{H}_2(\xi)&=&\lambda_0\Tr\{\mathcal{J}(\xi)\}=\lambda_0R(\xi,e_i,e_i,\xi),\\
\mathcal{H}_3(\xi)&=&\lambda_1\Tr\{\nabla_\xi\mathcal{J}(\xi)\}=\lambda_1R(\xi,e_i,e_i,\xi;\xi),\\
\mathcal{H}_4(\xi)&=&\lambda_2\Tr\{\nabla_\xi\nabla_\xi\mathcal{J}(\xi)\}+
\lambda_3\Tr\{\mathcal{J}(\xi)^2\}+\lambda_4\Tr\{\mathcal{J}(\xi)\}^2\\
&=&\lambda_2R(\xi,e_i,e_i,\xi;\xi\xi)+\lambda_3R(\xi,e_i,e_j,\xi) R(\xi,e_j,{e_i},\xi)\\
&&+\lambda_4R(\xi,e_i,e_i,\xi)R(\xi,e_j,e_j,\xi)
\end{array}\end{equation}
and so forth. To determine $\mathcal{H}_2$, $\mathcal{H}_3$, and $\mathcal{H}_4$,
we must determine the universal constants $\lambda_i$. The expressions for the $\mathcal{H}_i$ for
$i=5,6,7,8$ are, of course, more complicated.
\subsection{Universal examples}
Gray~\cite{Gr74} made a direct computation using curvature to determine $\mathcal{H}_n$ for $2\le n\le 6$.
Instead of using this approach, we use the method of universal examples. We suppose
given smooth functions $\Theta_i(r)$ for $1\le i\le m-1$ where the $\Theta_i$ are real analytic and $\Theta_i(0)=0$.
We consider the metric
$$
ds^2=dr^2+\sum_{i=1}^{m-1}\Theta_i^2d\theta_i^2\,.
$$
What we are doing, of course, is to take a parallel frame along a geodesic and
define the coordinates $\theta_i$ using the relevant geodesic spray; the origin is, as
always, singular in polar coordinates which is why we took the warping functions $\Theta_i(0)=0$.

We then have $\Theta_{\mathcal{M}}=\Theta_1\dots\Theta_{m-1}$.
For example, if we were to take all the functions $\Theta_i(r)=\sin(r)$, we
would obtain a representation of the metric on $S^m$ in geodesic polar coordinates
and if we were to take all the $\Theta_i(r)=\sinh(r)$, we
would obtain a representation of the metric on $\mathbb{H}^m$ in geodesic polar coordinates.
\subsection{Christoffel symbols} The coordinates decouple; to simplify the discussion, we take $m=2$. Let $f(r)=\Theta(r)^2$. Let $f_r:=\partial_rf$, $f_{rr}=\partial_r^2f$, and so forth.
We use the Koszul formula to compute:
$$\begin{array}{llll}
\Gamma_{rrr}=0,&\Gamma_{rr\theta}=0,&\Gamma_{rr}{}^r=0,&\Gamma_{rr}{}^\theta=0,\\[0.05in]
\Gamma_{r\theta r}=0,&\Gamma_{r\theta\theta}=\frac12{ f_{r}},&\Gamma_{r\theta}{}^r=0,
&\Gamma_{r\theta}{}^\theta=\frac12{ f_{r}}f^{-1},\\[0.05in]
\Gamma_{\theta\theta r}=-\frac12f_{r},&
\Gamma_{\theta\theta\theta}=\frac12f_{\theta},&\Gamma_{\theta\theta}{}^r=-\frac12f_{r},
&\Gamma_{\theta\theta}{}^\theta=\frac12{ f_{ \theta}} f^{-1}=0.
\end{array}$$
Thus we have that
\begin{eqnarray*}
&&\nabla_\theta \nabla_r \partial_r=0,\\
&& \nabla_r\nabla_\theta\partial_r={ \nabla_r}\{\Gamma_{\theta r}{}^\theta\partial_\theta\}=
{ \textstyle\{\frac12f_{rr}f^{-1}-\frac12f_{r}f_{r}f^{-2}+\frac14f_{r}f_{r}f^{-2}\}}\\
&&R(\partial_\theta,\partial_r,\partial_r,\partial_\theta)={ \textstyle-\frac12f_{rr}+\frac14f_{r}f_{r}f^{-2}\,.}
\end{eqnarray*}

More generally, we must subscript. Set
\begin{equation}\label{E2.b}\begin{array}{l}
\Theta_i(r)=r(1+b_{2,i}r^2+b_{3,i}r^3+\dots),\qquad f_i(r)=\Theta_i(r)^2,\\[0.05in]
 \Theta:=\prod_i\Theta_i(r),\qquad f_{i,r}:=\partial_rf_i,\qquad f_{i,rr}:=\partial_r^2f_i\,.
\end{array}\end{equation}
We have
\begin{equation}\label{E2.c}\begin{array}{l}
R(\partial_{\theta_i},\partial_r,\partial_r,\partial_{\theta_i} ) =
\textstyle-\frac12f_{i,rr}+\frac14{ f_{i,r}}f_{i,r}f_i^{-2},\\[0.05in]
R(\partial_{\theta_i},\partial_r,\partial_r,\partial_{\theta_j})=0\quad\text{for}\quad i\ne j,\\[0.05in]
 \nabla R(\partial_{\theta_i},\partial_r,\partial_r,\partial_{\theta_i};\partial_r)
=\partial_rR(\partial_{\theta_i},\partial_r,\partial_r,\partial_{\theta_i})-2\Gamma_{r\theta}{}^\theta R(\partial_{\theta_i},\partial_r,\partial_r,\partial_{\theta_i})
\end{array}\end{equation}
and so forth.
Since the $\{\partial_{\theta_i}\}$ do not form an orthonormal basis, we normalize by multiplying by $f_i^{-1}$ to
define
$$
\Lambda_{i,k}=f_i^{-1}\nabla^kR(\partial_{\theta_i},\partial_r,\partial_r,\partial_{\theta_i};\partial_r,\dots,\partial_r)\,.
$$
Thus, the eigenvalues of the operator
$\mathcal{J}_k(\partial_r)$ are $\{0,\Lambda_{1,k},\dots,\Lambda_{m-1,k}\}$ so
$$
\Tr\{\mathcal{J}_k({\partial_r})\}=\sum_{i=1}^{m-1}\Lambda_{i,k}\,.
$$

This procedure is ideally suited for implementing in a system of computer algebra.
Using this formalism gives rise to a number of equations when can then be solved to determine the unknown universal
coefficients. We took $m=5$ as this sufficed to determine $\mathcal{H}_i$ for $2\le i\le 8$;
 a more complicated
formalism would be needed to study $\mathcal{H}_9$ as these examples
can not distinguish between $\Tr\{\mathcal{J}\mathcal{J}_1\mathcal{J}_2\}$
and $\Tr\{\mathcal{J}\mathcal{J}_2\mathcal{J}_1\}$.
{One then uses Equation~(\ref{E2.b}) and Equation~(\ref{E2.c}) to compute:}
$$\begin{array}{l}
 \Tr\{\mathcal{J}(\partial_r)\}=-6(b_{2,1}+b_{2,2}+b_{2,3}+b_{2,4}),\\[0.05in]
 \mathcal{H}_2=b_{2,1}+b_{2,2}+b_{2,3}+b_{2,4}\,.
\end{array}$$
It then follows that the universal constant $\lambda_0$ of Equation~(\ref{E2.a}) is $-\frac16$ and we have
determined $\mathcal{H}_2$. Of course, as one goes further and further up in the expansion, the algebra
becomes more and more difficult and the use of a system of computer algebra becomes essential. Theorem~\ref{T1.1}
now follows.\qed

\begin{remark}\rm
If $\mathcal{M}$ is the round sphere of dimension $m$,
then
$$
\Theta=\sin(r)^{m-1},\quad\nabla R=0,\quad\text{and}\quad\Tr\{\mathcal{J}(\xi)^\ell\}=m-1\quad\text{for any}\quad\ell\,.
$$
We verified that the
formulas of Theorem~\ref{T1.1} are valid for this situation.
Similarly if $\mathcal{M}=\mathbb{CP}^k$, then $\Theta=\sin(r)^{2k-1}\cos(r)$, $\nabla R=0$, and $\Tr\{\mathcal{J}(\xi)^\ell\}=4^\ell+m-2$
and again the formulas are valid. A similar check holds for $\mathcal{M}=\mathbb{QP}^k$. This provides an independent
check that the coefficients of the terms
$\Tr\{\mathcal{J}(\xi)^{\ell_1}\}\dots\Tr\{\mathcal{J}(\xi)^{\ell_\nu}\}$ for $2\ell_1+\dots+2\ell_\nu=2k$ in $\mathcal{H}_{2k}$ are in fact correct.
\end{remark}
\begin{remark}\rm In our basic algorithm, we took $\Theta_i=r(1+b_{2,i}r^2+\dots)$. One can in fact take
$\Theta_i=r({b_{0,i}}+b_{2,i}r^2+\dots)$ and then one gets
$$
\Theta(r)=\prod_i b_{0,i}(1+\mathcal{H}_2r^2+\dots)\,.
$$
By taking $b_{0,i}$ to be $+1$ or $\sqrt{-1}$ and by taking the remaining $b_{2,i}$
to be real or pure imaginary, as appropriate, one can create a metric which is no longer positive definite but where $\xi$ is spacelike
and deduce that the relevant expansion holds, where the volume element is normalized appropriately;
changing the sign of the
metric then permits one to interchange spacelike with timelike vectors. We omit details in the interests of brevity.
\end{remark}

\section{The density function for a harmonic space}\label{S3}
In this section, we will use the formulas of Theorem~\ref{T1.1} to establish Theorem~\ref{T1.4}.
Let $\xi$ be a unit vector.
Since $\mathcal{M}$ is a harmonic space, $\mathcal{H}_2$ is constant, and hence
$\mathcal{M}$ is Einstein. Let $\sigma_\xi(r)$ be a geodesic with $\dot\sigma_\xi(0)=\xi$. The covariant derivative
of the trace is the trace of the covariant derivative.
\begin{equation}\label{E3.a}
\Tr\{\mathcal{J}_k(\xi)\}=\partial_r^k\Tr\{\mathcal{J}({\dot\sigma_\xi(r)})\}=0\quad\text{ for }\quad k\ge1\,.
\end{equation}
We impose the relations of Equation~(\ref{E3.a}) to show both that $\mathcal{H}_3=0$ and
that the formula of Assertion~2 for $\mathcal{H}_4$ holds. It now follows that
$\Tr\{\mathcal{J}(\xi)^2\}$ is constant and hence $\mathcal{M}$ is 2-stein.
We use Leibnitz's rule to differentiate $\Tr\{{ \mathcal{J}_{k}}(\xi)^2\}$
and obtain
\begin{equation}\label{E3.b}
\sum_{\mu+\nu=k}\frac{k!}{\mu!\nu!}\Tr\{\mathcal{J}_\mu(\xi)\mathcal{J}_\nu(\xi)\}=0\quad\text{ for }\quad k\ge1\,.
\end{equation}
We impose these relations; this shows $\mathcal{H}_5$ vanishes  and derives the formula given for $\mathcal{H}_6$.
Since $\mathcal{H}_6$ is constant and $\mathcal{M}$ is 2-stein,{
$$
-\frac1{2835}\Tr\{J(\xi)^3\}+\frac1{10080}\Tr\{\mathcal{J}_1(\xi)^2\}
=-\frac{32}{90720}\Tr\{J(\xi)^3\}+\frac9{90720}\Tr\{\mathcal{J}_1(\xi)^2\}
$$
is constant. We differentiate
$\Tr\{-32\mathcal{J}(\xi)^3+9\mathcal{J}_1^2(\xi)\}=\text{constant}$ to obtain
\begin{eqnarray*}
&&\Tr\{ -16\mathcal{J}(\xi)^2{\mathcal{J}_1(\xi)}+3\mathcal{J}_1(\xi)\mathcal{J}_2(\xi)\}=0,\\
&&\Tr\{{ -32\mathcal{J}(\xi)\mathcal{J}_1(\xi)^2- 16 \mathcal{J}(\xi)^2\mathcal{J}_2(\xi)+3\mathcal{J}_2(\xi)^2}
    +3\mathcal{J}_1(\xi)\mathcal{J}_3(\xi)\}=0\,.
\end{eqnarray*}
We impose these relations to see $\mathcal{H}_7=0$ and to determine $\mathcal{H}_8$.
This completes the proof of Theorem~\ref{T1.4}.\qed}

\section{Estimates on the eigenvalues}\label{S4}
In Section~\ref{S4.1}, we establish Lemma~\ref{L1.7}, in Section~\ref{S4.2}, we demonstrate Theorem~\ref{T1.8},
and in Section~\ref{S4.3} we prove Theorem~\ref{T1.9}.
In Section~\ref{S4.4} we indicate why the structure of the density function, the Jacobi operator, and the structure of the Jacobi
vector fields are so closely interrelated.
Throughout this section, we shall assume that $\mathcal{M}$ is a harmonic space
satisfying $\Theta_{\mathcal{M}}=\Theta_\varepsilon(m,k)$ where $\Theta_\varepsilon(m,k)$ is as defined
in Equation~(\ref{E1.b}).
Order the eigenvalues of the reduced Jacobi operator $\tilde{\mathcal{J}}_{\mathcal{M}}(\xi)$ so that
$\varepsilon\lambda_1\ge\dots\ge\varepsilon_{m-1}\lambda_{m-1}$. The following identities will be crucial.
\begin{lemma}\label{L4.1}
Assume that $\Theta_{\mathcal{M}}=\Theta_\varepsilon(m,k)$.  Then
$$
\sum_{i=1}^{m-1}\varepsilon\lambda_i=4k+(m-k-1)\text{ and }
\sum_{i=1}^{m-1}\lambda_i^2=16k+(m-k-1)\,.
$$
\end{lemma}

\begin{remark}\rm
These are exactly the values that would pertain if $\mathcal{M}$ was a rank 1 symmetric space.
\end{remark}

\begin{proof}
We expand
\medbreak\qquad$\sin(r)=r-\frac16r^3+\frac1{120}r^5-\frac1{5040}r^7+O(r^9)$,
\smallbreak\qquad$\cos(r)=1-\frac12r^2+\frac1{24}r^4-\frac1{720}r^6+O(r^8)$,
\smallbreak\qquad$\sinh(r)=r+\frac16r^3+\frac1{120}r^5+\frac1{5040}r^7+O(r^9)$,
\smallbreak\qquad$\cosh(r)=1+\frac12r^2+\frac1{24}r^4+\frac1{720}r^6+O(r^8)$.
\medbreak\noindent
It is then a direct computation to see:
\begin{equation}\begin{array}{l}\label{E4.a}
\Theta_\varepsilon(m,k)=r^{m-1}\left(1+\varepsilon\frac{-3 k-m+1}{6} r^2 +\frac{45 k^2+30 k m-60 k+5 m^2-12 m+7}{360} r^4\right.\\[0.05in]
\quad+\varepsilon\frac{-945 k^3-945 k^2 m+2835 k^2-315 k m^2+1386 k m-2079 k-35 m^3+147 m^2-205 m+93}{45360}r^6\\[0.05in]
\quad+O(r^8)\left.\vphantom{A^{A^A}}\right)\,.
\end{array}\end{equation}
We have $\Tr\{\mathcal{J}(\xi)\}=\sum_i\lambda_i$ and $\Tr\{\mathcal{J}^2(\xi)\}=\sum_i\lambda_i^2$. By Theorem~\ref{T1.4}
and Equation~(\ref{E4.a}),
\begin{eqnarray*}
&&\mathcal{H}_2=-\frac16\Tr\{\mathcal{J}(\xi)\}=\varepsilon\frac{-3 k-m+1}{6},\\
&&\mathcal{H}_4=\frac1{72}\Tr\{\mathcal{J}\}^2-\frac1{180}\Tr\{\mathcal{J}^2\}=\frac{45 k^2+30 k m-60 k+5 m^2-12 m+7}{360}\,.
\end{eqnarray*}
We solve these relations for $\Tr\{\mathcal{J}\}$ and $\Tr\{\mathcal{J}^2\}$ to establish the Lemma.
\end{proof}

\subsection{The proof of Lemma~\ref{L1.7}}\label{S4.1}

As $g_{ij}(\xi)=\delta_{ij}+O(\|\xi\|^2)$ in geodesic normal coordinates,
$\exp_P^*(\operatorname{dvol}_{\mathcal{M}})=\operatorname{dvol}_{\mathcal{M}}$ at the origin.
Consequently,
$$
\operatorname{dvol}_{\mathcal{M}}=r^{\dim(M)-1}drd\theta\quad\text{so}\quad
\Theta_{\mathcal{M}}(P;\xi)=\|\xi\|^{\dim(M)-1}+O(\|\xi\|^{\dim(M)})\,.
$$
By Equation~(\ref{E4.a}), $\Theta_\varepsilon(m,k)(r)=r^{m-1}+O(r^m)$.
Thus
$$
\Theta_{\mathcal{M}}=\Theta_\varepsilon(m,k)\quad\text{implies}\quad\dim(M)=m\,.
$$
Set $\alpha:=\lambda_1+\dots+\lambda_{m-1}$. Express $\lambda_i=\frac\alpha{m-1}+\delta_i$.
We then have $\delta_1+\dots+\delta_{m-1}=0$ so
\begin{equation}\label{E4.b}\begin{array}{l}
\displaystyle\sum_{i=1}^{m-1}\lambda_i^2=\sum_{i=1}^{m-1}\left(\frac\alpha{m-1}+\delta_i\right)^2\\[0.05in]
\displaystyle\phantom{...}\qquad=\frac{(m-1)\alpha^2}{(m-1)^2}+\frac{2\alpha}{m-1}\sum_{i=1}^{m-1}\delta_i+\sum_{i=1}^{m-1}\delta_i^2\ge\frac{\alpha^2}{m-1}\,.
\end{array}\end{equation}
Consequently, we may rewrite Equation~(\ref{E4.b}) as the estimate:
\begin{equation}\label{E4.c}
\sum_{i=1}^{m-1}\lambda_i^2\ge\frac1{m-1}\left(\sum_{i=1}^{m-1}\lambda_i\right)^2\,.
\end{equation}
We apply Lemma~\ref{L4.1} and Equation~(\ref{E4.c}) to estimate
\begin{equation}\label{E4.d}
16k+m-k-1\ge\frac{(4k+m-k-1)^2}{m-1}\,.
\end{equation}
We set
$f(t):=(4t+m-t-1)^2-(m-1)(16t+m-t-1)=9t(t-(m-1))$.
Since $f(k)\le0$, $0\le k\le m-1$. This completes the proof of Lemma~\ref{L1.7}~(1).

We now examine the extremal cases. If $k=0$ or if $k=m-1$, then we have equality in Equation~(\ref{E4.d}). This means we have
equality in Equation~(\ref{E4.c}) and hence we have equality in Equation~(\ref{E4.b}). This implies that all the $\delta_i$ vanish and hence
$\lambda_i=\frac\alpha{m-1}$ so all the eigenvalues of the reduced Jacobi operator are constant. This implies
that $\mathcal{M}$ has constant sectional curvature; the sectional curvature is $\varepsilon$ if $k=0$ and $4\varepsilon$ if $k=m-1$ by
Lemma~\ref{L4.1}.
\qed
\subsection{The proof of Theorem~\ref{T1.8}}\label{S4.2}
Express
\begin{equation}\label{E4.e}
\varepsilon\lambda_i=\left\{\begin{array}{ll}4+\delta_i&\text{if }1\le i\le k\\1+\delta_i&\text{if }k+1\le i\le m-1\end{array}\right\}\,.
\end{equation}
By Lemma~\ref{L4.1} and Equation~(\ref{E4.e}), we have
\begin{equation}\label{E4.f}
\sum_{i=1}^{m-1}\delta_i=0\,.
\end{equation}
We now compute:
\begin{equation}\label{E4.g}\begin{array}{l}
\displaystyle16k+(m-k-1)=\sum_{i=1}^{m-1}\lambda_i^2=\sum_{i=1}^k(4+\delta_i)^2+\sum_{i=k+1}^{m-1}(1+\delta_i)^2\\[0.05in]
\displaystyle\qquad=16k+(m-k-1)+8\sum_{i=1}^k\delta_i+2\sum_{i=k+1}^{m-1}\delta_i+\sum_{i=1}^{m-1}\delta_i^2
\end{array}\end{equation}
Using Equation~(\ref{E4.f}) and Equation~(\ref{E4.g}), we see
\begin{equation}\label{E4.h}
6\sum_{i=1}^k\delta_i+\sum_{i=1}^{m-1}\delta_i^2=0\quad\text{so}\quad\sum_{i=1}^k\delta_i\le0\,.
\end{equation}
This implies
\begin{equation}\label{E4.i}
\sum_{i=1}^k \varepsilon\lambda_i=\sum_{i=1}^k(4+\delta_i)\le 4k\,.
\end{equation}
We complete the proof of Theorem~\ref{T1.8}~(1)
by using Lemma~\ref{L4.1} and Equation~(\ref{E4.i}) to estimate
\begin{equation}\begin{array}{l}\label{E4.j}
\displaystyle\sum_{i=k+1}^{m-1}\varepsilon\lambda_i=\sum_{i=1}^{m-1}\varepsilon\lambda_i-\sum_{i=1}^k\varepsilon\lambda_i=4k+m-1-k-\sum_{i=1}^k\varepsilon\lambda_i\\
\qquad\qquad\vphantom{.......}\ge4k+m-1-k-4k=m-1-k\,.
\end{array}\end{equation}

We suppose either (and hence both) of the inequalities of Equation~(\ref{E4.i}) or Equation~(\ref{E4.j}) are in fact equalities. This implies $\sum_{i=1}^k\delta_i=0$ and hence by Equation~(\ref{E4.h}),
all the $\delta_i=0$. Consequently
$$
\varepsilon\lambda_i=\left\{\begin{array}{ll}4&\text{if }1\le i\le k\\1&\text{if }k+1\le i\le m-1\end{array}\right\}\,.
$$
Since this is assumed true for every unit tangent vector, we conclude $\mathcal{M}$ is Osserman. We may now
compute $\operatorname{Tr}(\mathcal{J}_{\mathcal{M}}^i)$ for $i=1,2,3$.
We equate the $6^{\operatorname{th}}$ order term in the Taylor series for $\Theta_\pm(m,k)(r)$ with $\mathcal{H}_6$,
to see $\frac1{10080}\operatorname{Tr}((\nabla_\xi\mathcal{J}_{\mathcal{M}}(\xi))^2)=0$.
Since
$\nabla_\xi \mathcal{J}_{\mathcal{M}}(\xi)$ is self-adjoint, this implies $\nabla_\xi \mathcal{J}_{\mathcal{M}}(\xi)=0$.
Since this assumed true for every $\xi$, we obtain that the matrix of $\mathcal{J}_{\mathcal{M}}(\dot\sigma)$ is constant along any geodesic
relative to a parallel frame. This implies that the geodesic involution is necessarily an isometry and hence $\mathcal{M}$
is locally symmetric (see~\cite{B78} Proposition 2.35).

We are working locally so we may assume $M$ is a small geodesic ball. Since $M$ is
symmetric, it follows that $M$ embeds isometrically in a globally symmetric simply connected
manifold, see Helgason~\cite{H60} Theorem 5.1.
Heber~\cite{H06} showed that a simply connected homogeneous harmonic space is
flat, or is a rank-one symmetric space, or is a non-symmetric Damek-Ricci space. Since $\mathcal{M}$
is not flat and is symmetric, this completes the proof.\qed
\subsection{Bounding the spectrum from above}\label{S4.3}
We can use a similar argument to bound the spectrum. We
drop our convention that the eigenvalues are ordered according to size. We consider the problem of minimizing/maximizing
the function $\varepsilon\lambda_1$ subject to the constraints
$$
\varepsilon\lambda_1+\dots+\varepsilon\lambda_{m-1}=A\quad\text{and}\quad\lambda_1^2+\dots+\lambda_{m-1}^2=B\,.
$$
Using  the method of Lagrange multipliers then yields
$$
(1,0,\dots,0)=\mu_1(1,1,\dots,1)+\mu_2(2\lambda_1,\dots,2\lambda_{m-1})\,.
$$
Since $\mu_2\ne0$, we obtain $\lambda_2=\dots=\lambda_{m-1}$. The first constraint then yields
$$
\varepsilon\lambda_i=\frac{A-\varepsilon\lambda_1}{m-2}\quad\text{for}\quad2\le i\le m-1\,.
$$
Substituting this into the second constraint then yields $\lambda_1^2+\frac{(A-\varepsilon\lambda_1)^2}{m-2}=0$. Solving the
resulting equation with the values for $A$ and $B$ given by Lemma~\ref{L4.1} then yields two roots; the larger provides an upper bound
on an eigenvalue
$$
\varepsilon\lambda_i\le(m-1)^{-1}\left\{-1+3k+m+(3m-6)\sqrt{\frac{k(m-1-k)}{m-2}}\right\}\,.
$$
and the smaller root provides a lower bound on an eigenvalue
$$
\varepsilon\lambda_i\ge(m-1)^{-1}\left\{-1+3k+m-(3m-6)\sqrt{\frac{k(m-1-k)}{m-2}}\right\}\,.
$$
For example, if $m=8$ and $k=3$, we obtain $\varepsilon\lambda_i\in[-1.35084,5.92227]$. This completes the proof of
Theorem~\ref{T1.9}.

\subsection{The index form}\label{S4.4}
We digress and follow the discussion in \cite{GHL,L}.
The volume density function is intimately related to the structure of the Jacobi vector fields.
Let $Y$ be a vector field which
is defined along a unit speed geodesic $\sigma:[0,T]\rightarrow M$. We say that $Y$ is a {\it Jacobi vector field} if it satisfies
the identity $\nabla_{\dot\sigma}\nabla_{\dot\sigma}Y+\mathcal{J}_{\mathcal{M}}(\dot\sigma)Y=0$. Jacobi vector fields
arise by considering geodesic sprays. Let $\sigma_s(t)$ be a 1-parameter family of geodesics. We then have
that $\partial_s$ is a Jacobi vector field along any geodesic $\sigma_s(t)$. The {\it index form} is the symmetric bilinear form defined by setting
$$
I_\sigma(V,W):=\int_{t=0}^T\left\{g(\nabla_{\dot\sigma}V,\nabla_{\dot\sigma}W)-R(V,\dot\sigma,\dot\sigma,W)\right\}dt\,.
$$
The following result relates the index form to the volume density function and the mean curvature of the geodesic spheres:
\begin{lemma}
Let $\sigma:[0,T]\rightarrow M$ be a unit speed geodesic in a Riemannian manifold $\mathcal{M}$.
Assume $T$ is small so there are no conjugate points along $\sigma$.
Let $\{Y_1,\dots,Y_{m-1}\}$
be perpendicular Jacobi vector fields along $\sigma$ with $\{\dot\sigma(T),Y_1(T),\dots,Y_{m-1}(T)\}$ an orthonormal
basis for $T_{\sigma(T)}$ and with $Y_i(0)=0$ for $1\le i\le m-1$.
Let $\mathfrak{m}$ be the mean curvature of the geodesic sphere about $\sigma(0)$ passing thru $\sigma(T)$. Then
$$
\mathfrak{m}=\frac{\dot\Theta_P(\sigma(T))}{\Theta_P(\sigma(T))}=\sum_{i=1}^{m-1}I_\sigma(Y_i,Y_i)\,.
$$
\end{lemma}
Although this observation played no role in our development, it formed the starting point for the analysis of
Ramachandran and Ranjan~\cite{RR97} and one can use this result to give an alternate derivation of Theorem~\ref{T1.8}
based upon their methods; we
omit details in the interests of brevity.

\section*{Research Support} This work was supported by the National Research Foundation of Korea (NRF) grant
funded by the Korea government(MSIT) (NRF-2019R1A2C1083957) and by Project MTM2016-75897-P (AEI/FEDER, UE). The authors acknowledge with gratitude useful
conversations with Professor Garcia-Rio concerning the paper.


\begin{thebibliography}{99}
\bibitem{AD}J-P Anker, E. Damek, C.  Yacoub, ``Spherical analysis on harmonic
AN groups", {\it Ann. Scuola Norm. Sup. Pisa Cl. Sci. \bf 23} (1996), 643--679 (1997).
\bibitem{BTV78}J. Berndt, F. Tricerri and L. Vanhecke, ``Generalized Heisenberg Groups and Damek Ricci Harmonic Spaces", Lecture Notes in Math. 1598, Springer-Verlag, Berlin 1995.
\bibitem{B78}A.L. Besse, ``Manifolds all of whose geodesics are closed", Ergebnisse der Mathematik
und ihrer Grenzgebiete 93, Springer-Verlag, Berlin 1978.
\bibitem{BG95}Besson, G. Courtois, S. Gallot, ``Entropies et rigidities des espaces localement
symtriques de courbure strictement negative", {\it Geometric and functional analysis \bf5}, (1995).
\bibitem{CGW82} P. Carpenter,  A. Gray, and T. J. Willmore,
``The curvature of Einstein symmetric spaces", {\it Quart. J. Math. Oxford, \bf33} (1982), 45--64.
\bibitem{C88} C. Chi, ``A curvature characterization of certain locally rank-one symmetric spaces",
{\it  J. Differential Geom. \bf 28} (1988), 187--202.
\bibitem{CR40}E. T. Copson and H. S. Ruse, ``Harmonic Riemannian spaces", {\it Proc. Roy. Soc. Edinburgh \bf60} (1940), 117--133.
\bibitem{CH16} Csik\'os and Horv\'ath, ``Harmonic manifolds and the volume of tubes about curves",
{\bf J. London Math. Soc. \bf94} (2016), 141--160.
\bibitem{CH18} B. Csik\`os and M. Horv\'ath, ``Harmonic manifolds and tubes", {\it J. Geom. Anal. \bf 2018}, 3458--3476.
\bibitem{DR57}E. Damek and F. Ricci,  ``A class of nonsymmetric harmonic Riemannian spaces",
{\it Bull. Amer. Math. Soc. \bf27} (1992), 139--142.
\bibitem{EPS}Y. Euh, J. H. Park, and K. Sekigawa, ``Characteristic function of Cayley projective plane as a harmonic manifold",
{\it Hokkaido Math. J. \bf47} (2018), 191-203.
\bibitem{GHL} S. Gallot, D. Hulin, J. Lafontaine, {\bf Riemannian geometry, Third edition}, Universitext, Springer-Verlag, Berlin, 2004, xvi+322 pp.
\bibitem{G95} P. Gilkey, {\bf Invariance theory,
the heat equation, and the Atiyah-Singer Index Theorem} $2^{\operatorname{nd}}$ ed, CRC Press (Boca Raton), 1995.
\bibitem{Gr74}A. Gray, ``The volume of a small geodesic ball of a Riemannian manifold", {\it Michigan Math. J. \bf 20} (1973), 329--344.
\bibitem{H06}J. Heber, ``On harmonic and asymptotic harmonic homogeneous spaces", {\it Geom. Funct. Anal. \bf16}
(2006), 869--890.
\bibitem{H60}S. Helgason, {\bf Differential Geometry and Symmetric Spaces}, Pure and Applied mathematics {\bf XII}, Academic Press (New York--London) (1962).
\bibitem{KP} S. Kim and J.H. Park, ``Characteristic function of a symmetric Damek Ricci space",
{\it Bull. Korean Math Soc. \bf 54} (2017), 911--916.
\bibitem{L} J. M. Lee, {\bf Riemannian manifolds}, Graduate Texts in Mathematics, 176, Springer, New York, Inc.,1997.
\bibitem{L57}A. J. Ledger, ``Symmetric harmonic spaces", {\it J. London Math. Soc. \bf 32} (1957), 53--56.
\bibitem{L44}A. Lichnerowicz, ``Sur les espaces riemanniens completement harmoniques",
{\it Bull. Soc. Math. France \bf72} (1944), 146--168.
\bibitem{nikox} Y. Nikolayevsky, ``Harmonic homogeneous manifolds of nonpositive curvature",\newline  arXiv:math/0407024.
\bibitem{niko3} Y. Nikolayevsky,  ``Osserman manifolds of dimension 8", {\it Manuscripta Math. \bf 115} (2004), 31--53.
\bibitem{niko4} Y. Nikolayevsky, ``Osserman conjecture in dimension {$n\neq 8,16$}"
{\it Math. Ann. \bf 331} (2005), 505--522.
\bibitem{N05}Y. Nikolayevsky, ``Two theorems on harmonic manifolds", {\it Comment. Math. Helv. \bf80} (2005), 29--50.
\bibitem{oss} R. Osserman, ``Curvature in the eighties", {\it Amer. Math. Monthly \bf 97} (1990), 731--756.
 \bibitem{RR97}K. Ramachandran and A. Ranjan, ``Harmonic manifolds with some specific volume
densities", {\it Proc. Indian Acad. Sci. Math. Sci. \bf107} (1997),
251--261.
\bibitem{RS2002} A. Ranjan and H. Shah, ``Harmonic manifolds with minimal horospheres", {\it J. Geom.
Anal. \bf12} (2002), 683--694.
\bibitem{S90}Z. I. Szab\'o, ``The Lichnerowicz conjecture on harmonic manifolds", {\it J. Differential Geom. \bf31} (1990), 1--28.
\bibitem{S93}Z. I. Szab\'o, ``Spectral theory for operator families on Riemannian manifolds", {\it Proc. Sympos. Pure Math. \bf54}, Part 3, Amer. Math. Soc., Providence, RI, 1993.
\bibitem{V}L. Vanhecke, ``A note on harmonic spaces",
{\it Bull. London Math Soc. \bf 13} (1981), 545--546.
\bibitem{W49}A. G. Walker, ``On Lichnerowicz's conjecture for harmonic 4-spaces", {\it J. London Math. Soc. \bf24} (1949), 21--28.
\bibitem{W75} Y. Watanabe, ``On the characteristic function of harmonic K\"ahlerian spaces'', {\it Tohoku Math. J. \bf27} (1975), 13--24.
\bibitem{W76} Y. Watanabe,
``On the characteristic functions of quaterion K\"ahlerian spaces", {\it Kodai Math. Sem. Rep. \bf 27} (1976),
410--420.

\end{thebibliography}
\end{document}